\documentclass[preprints,article,accept,moreauthors,pdftex]{my-mdpi} 

\firstpage{1} 
\pubvolume{10}
\issuenum{4}
\articlenumber{314}
\pubyear{2021}
\copyrightyear{2021}
\externaleditor{Academic Editor: Gabriella Bretti} 
\datereceived{23 October 2021} 
\daterevised{19 November 2021} 
\dateaccepted{21 November 2021} 
\datepublished{23 November 2021} 
\hreflink{https://doi.org/10.3390/\newline axioms10040314} 
\doinum{10.3390/axioms10040314}

\pdfoutput=1

\Title{A Discrete-Time Compartmental Epidemiological Model for COVID-19 with a Case Study for Portugal}

\TitleCitation{A~Discrete-Time Compartmental Epidemiological Model for COVID-19 with a Case Study for Portugal}

\Author{Sandra Vaz $^{1,\dagger}$\orcidA{} 
and Delfim F. M. Torres $^{2,}$*$^{,\dagger}$\orcidB{}}

\AuthorNames{Sandra Vaz and Delfim F. M. Torres}

\AuthorCitation{Vaz, S.; Torres, D.F.M.}

\address{$^{1}$ \quad Center of Mathematics and Applications (CMA-UBI),
Department of Mathematics, University of Beira Interior, 6201-001 Covilh\~{a}, Portugal; svaz@ubi.pt\\
$^{2}$ \quad Center for Research and Development in Mathematics and Applications (CIDMA),
Department of Mathematics, University of Aveiro, 3810-193 Aveiro, Portugal}

\corres{Correspondence: delfim@ua.pt; Tel.: +351-234-370-668}

\firstnote{These authors contributed equally to this work.}

\abstract{In [Ecological Complexity 44 (2020) Art. 100885, DOI: 10.1016/j.ecocom.2020.100885]
a continuous-time compartmental mathematical model for the spread of the
Coronavirus disease 2019 (COVID-19) is presented with Portugal as case study, 
from 2 March to 4 May 2020, and the local stability of the 
Disease Free Equilibrium (DFE) is analysed. Here, we propose an analogous discrete-time model 
and, using a suitable Lyapunov function, we prove the global stability of the DFE point. 
Using COVID-19 real data, we show, through numerical simulations, 
the consistence of the obtained theoretical results.}

\keyword{mathematical modelling; epidemiology; COVID-19 pandemic; qualitative theory; \linebreak
\mbox{Lyapunov functions}; numerical simulations with real data}

\MSC{65L12; 92D30}

\begin{document}


\section{Introduction}

The coronavirus belong to the family of \emph{Coronaviridae}. It is a virus that causes infection 
to humans, other mammals and birds. The infection usually affects the respiratory system,
 and the symptoms may vary from simple colds to pneumonia. To date, eight coronaviruses are known. 
The new coronavirus SARS-CoV-2 originates COVID-19, and was identified for the first time 
in December 2019 in Wuhan, China. Two other coronaviruses have caused outbreaks: SARS-CoV, 
in 2002--2003, and MERS-CoV in 2012 \citep{Covid19}. The COVID-19 pandemic in an ongoing pandemic, 
probably the most significant one in human history. The effects are so severe that it has been 
necessary to use quarantining to control the spread. It has forced humans to adapt 
to this new situation \cite{MyID:461}.

The use of mathematical compartmental models to study the spread and consequences  
of infectious diseases have been used successfully for a long time. 
The techniques to study compartmental models are vast: continuous models, fractional models, 
and discrete models, among others \cite{MyID:461}. The complexity of the models is determined 
by the effects and consequences of the disease. Several models have been proposed for 
COVID-19 spread; \mbox{see \cite{MyID:461,PST,ST,NT}}, among others. In this work, we are 
interested in discretizing the model presented in \cite{PST}, using the Mickens method
\cite{Mickens94,Mickens02}, which was applied successfully in several 
different contexts \cite{VD}.

\textls[-10]{In Section~\ref{sec2}, we recall the SAIQH 
(susceptible--asymptomatic--infectious--quarantined--}hospitalized) continuous-time 
mathematical model, the equilibrium points and the stability results of \cite{PST}. 
Section~\ref{sec3} is dedicated to the original results of our work: 
the new discrete-time model, the proof of its well-posedness, 
the computation of the equilibrium points, and the proof
of the stability of the disease-free equilibrium point. We end Section~\ref{sec3}
with numerical simulations showing the consistency of our results. 
Finally, Section~\ref{sec4} is dedicated to discussion
of the obtained results and some possible future work directions.  


\section{Preliminaries}
\label{sec2}

This section is dedicated to the presentation of the continuous-time model \cite{PST}
and its results, which establish the stability of the equilibrium points. 
For more information, \mbox{see \cite{PST}}.

The total living population under study at time $t \geq 0$ is denoted by $N(t)$ as follows: 
\begin{equation}
N(t)=S(t)+A(t)+I(t)+Q(t)+H(t)+\overline{H}(t) 
\end{equation}
where $S(t)$ represents the susceptible individuals; $A(t)$ the infected individuals 
without symptoms or with mild ones; $I(t)$  the infected individuals; $Q(t)$ the 
individuals in quarantine, that is, in isolation at home; $H(t)$ the hospitalized 
individuals; and, finally, $\overline{H}(t)$ the hospitalized individuals in intensive care units.
There is another class that is also considered, $D(t)$, that gives the cumulative 
number of deaths due to COVID-19 for all $t \geq 0$. 
Regarding the parameters, all of them are positive. 
The recruitment rate into the susceptible class is $\Lambda$; 
all individuals of all classes are subject to the natural 
death rate $\mu$ along all time $t \geq 0$ under study. 
Susceptible individuals may become infected with COVID-19 at the following rate: 
\begin{equation}
\lambda(t)=\frac{\beta (l_{A} A(t)+I(t)+l_{H}H(t))}{N(t)},
\end{equation}
where $\beta$ is the human-to-human transmission rate per unit of time (day) 
and $l_{A}$ and $l_{H}$ quantify the relative transmissibility of asymptomatic 
individuals and hospitalized individuals, respectively. The class $\overline{H}$ 
does not enter in $\lambda$ because the number of health care workers that become 
infected by SARS-CoV-2 in intensive care units is very low and can be neglected. 
A fraction $p \in [0,1]$ of the susceptible population is in quarantine at home, 
at rate $\phi$. Consequently, only a fraction $1-p \in [0,1]$ of susceptible individuals 
are assumed to be able to become infected. Since there is uncertainty about long-immunity 
after recovery, it is assumed that individuals of class $Q$ will become susceptible 
again at a rate $\omega$. It is also considered that only a fraction $m \in[0,1]$ 
of the quarantined individuals move from class $Q$ to $S$.  It means that $(m \times 100)\%$ 
of the quarantined individuals return to class $S$ at the end  of $\frac{1}{\omega}$ days. 
These assumptions are justified by the state of calamity that was immediately decreed 
by the government of Portugal to address the state of emergency, which was fully respected 
by the Portuguese population.  After $\frac{1}{\nu}$ days of infection, only 
a fraction $q \in [0,1]$ of infected individuals without (or with mild) symptoms 
have severe symptoms. Thus, $(q \times 100)\%$ of individuals  of compartment 
$A$ move to $I$ at rate $\nu$. A fraction $f_{1} \in [0,1]$ of infected individuals 
with severe symptoms are treated at home and the other fraction $(1-f_{1}) \in [0,1]$ 
are hospitalized, both at rate $\delta_{1}$. The model considers the three following scenarios 
for hospitalized individuals:
\begin{enumerate}
\item A fraction $f_{2} \in [0,1]$ of individuals in class $H$ can evolve 
to a state of severe health status, needing an invasive intervention, 
such as artificial respiration, so they need to move to intensive care, 
at rate $\delta_{2}$;

\item A fraction $f_{3} \in [0,1]$ of individuals in class $H$ die due to COVID-19, 
the disease related death rate associated with hospitalized individuals being $\alpha_{1}$;
 
\item A fraction $(1-f_{2}-f_{3}) \in [0,1]$ of individuals in class $H$ 
recover and, consequently, return home in quarantine/isolation 
at rate $\delta_{2}$.
\end{enumerate}

Regarding the hospitalized individuals in intensive care units, 
the model considers two possibilities as follows:
\begin{enumerate}
\item A fraction $(1-\kappa) \in [0,1]$ of individuals in class $\overline{H}$ 
recover and move to the class $H$ at rate $\eta$;

\item A fraction $\kappa \in [0,1]$ of individuals in class $\overline{H}$ 
die due to COVID-19, the disease-related death rate associated with 
hospitalized individuals in intensive care units being $\alpha_{2}$.
\end{enumerate}

Compiling all the previous assumptions, one has the following mathematical model:
\begin{equation}
\label{eq:model}
\begin{cases}
\dot{S}(t) = \Lambda + \omega m Q(t) - [ \lambda(t) (1-p) + \phi p + \mu] S(t),\\[0.2 cm]
\dot{A}(t) = \lambda(t)(1-p) S(t) - (q v +\mu) A(t), \\[0.2 cm]
\dot{I}(t) = q v A(t) - ( \delta_{1} + \mu)I(t),\\[0.2 cm]
\dot{Q}(t) = \phi p S(t) + \delta_{1} f_{1} I(t) + \delta_{2}(1-f_{2}-f_{3})H(t) 
- (\omega m + \mu) Q(t),\\[0.2 cm]
\dot{H}(t)=\delta_{1}(1-f_{1})I(t)+\eta(1-\kappa)\overline{H}(t)
-(\delta_{2}(1-f_{2}-f_{3})+\delta_{2}f_{2}+ \alpha_{1}f_{3}+\mu)H(t),\\[0.2 cm]
\dot{\overline{H}}(t)=\delta_{2}f_{2} H(t) 
- (\eta(1-k) + \alpha_{2} \kappa + \mu)\overline{H}(t), \\[0.2 cm]
\dot{D}(t)=\alpha_{1} f_{3} H(t)+ \alpha_{2} \kappa \overline{H}(t).
\end{cases}
\end{equation}

It can be seen in \cite{PST} that \eqref{eq:model} is equivalent to the following:
\begin{equation}
\label{eq:modelC}
\begin{cases}
\dot{S}(t) = \Lambda + \omega m Q(t) - [ \lambda(t) (1-p) + \phi p + \mu] S(t),\\[0.2 cm]
\dot{A}(t) = \lambda(t)(1-p) S(t) - (q v +\mu) A(t), \\[0.2 cm]
\dot{I}(t) = q v A(t) - ( \delta_{1} + \mu)I(t),\\[0.2 cm]
\dot{Q}(t) = \phi p S(t) + \delta_{1} f_{1} I(t) + \delta_{2}(1-f_{2}-f_{3})H(t) 
- (\omega m + \mu) Q(t),\\[0.2 cm]
\dot{H}(t)=\delta_{1}(1-f_{1})I(t)+\eta (1-\kappa)\overline{H}(t)-(\delta_{2}(1-f_{2}-f_{3})
+\delta_{2}f_{2}+ \alpha_{1}f_{3}+\mu)H(t), \\[0.2 cm]
\dot{\overline{H}}(t)=\delta_{2}f_{2} H(t) - (\eta(1-k) + \alpha_{2} \kappa + \mu)\overline{H}(t).
\end{cases}
\end{equation} 

Table~\ref{tab1} presents all parameters and initial conditions 
of mathematical model \eqref{eq:modelC}.

In \cite{PST}, it is shown that the biologically feasible region is given by the following:
\begin{equation}
\Omega=\left\{(S,A,I,Q,H,\overline{H}) 
\in ( \mathbb{R}_{0}^{+})^{6}: N \leq \frac{\Lambda}{\mu} \right\},
\end{equation}
which is positively invariant for \eqref{eq:modelC} for all 
non-negative initial conditions. To simplify the expressions 
in the computations, the following notation is used:
\begin{enumerate}
\item[(i)] $a_{0}:= q v + \mu$;
\item[(ii)] $a_{1}:= \delta_{1}+ \mu$;
\item[(iii)] $a_{2}:= m \omega + \mu$;
\item[(iv)] $a_{3}:= \delta_{2}(1-f_{2}-f_{3}) + \delta_{2}f_{2} + \alpha_{1} f_{3} + \mu$;
\item[(v)] $a_{4}:= \delta_{2}(1-f_{2}-f_{3})$;
\item[(vi)] $a_{5}:= p \phi + \mu$;
\item[(vii)] $a_{6}:= \delta_{1} (1 - f_{1})$;
\item[(viii)] $a_{7}:=\alpha_2 \kappa + \eta_{k} + \mu$;
\item[(ix)] $\eta_{k} := \eta(1-\kappa)$;
\item[(x)] $\chi :=a_{3} a_{7} - \delta_{2} \eta_{k} f_{2}$.
\end{enumerate}

\begin{specialtable}[H] 
\caption{Description of the parameters and initial conditions 
of model \eqref{eq:modelC}.\label{tab1}}
\begin{tabular}{cm{11.2cm}<{\raggedright}}
\toprule
\textbf{Parameter}	& \textbf{Description}\\
\midrule
$\Lambda$ & Recruitment Rate	  \\
$\mu$	  & Natural death rate	\\
$\beta$   & Human-to-human transmission rate \\
$l_{A}$   & Relative transmissibility of individuals in class $A$\\
$l_{H}$   & Relative transmissibility of individuals in class $H$\\
$\phi$    & Rate associated with movement from $S$ to $Q$\\
$\nu$     & Rate associated with movement from $A$ to $I$\\
$\delta_{1}$ & Rate associated with movement from $I$ to $Q/H$\\
$\delta_{2}$ & Rate associated with movement from $H$ to $Q/\overline{H}$\\
$\eta$       & Rate associated with movement from $\overline{H}$ to $H$\\
$\omega$     & Rate associated with movement from $Q$ to $S$\\
$\alpha_{1}$ & Disease-related death rate of class $H$\\
$\alpha_{2}$ & Disease-related death rate of class $\overline{H}$\\ 
$p$      & Fraction of susceptible individuals putted in quarantine\\
$q$      & Fraction of infected individuals with severe symptoms\\
$f_{1}$  & Fraction of infected individuals with severe symptoms in quarantine\\
$f_{2}$  & Fraction of hospitalized individuals transferred to $\overline{H}$\\
$f_{3}$  & Fraction of hospitalized individuals who die of COVID-19\\
$\kappa$ & Fraction of hospitalized individuals in intensive care units\\
& who die from COVID-19\\
$m$ & Fraction of individuals who moves from $Q$ to $S$\\
$S(0)=S_{0}$& Individuals in class $S$ at $t=0$\\
$A(0)=A_{0}$ & Individuals in class $A$ at $t=0$\\
$I(0)=I_{0}$ & Individuals in class $I$ at $t=0$\\
$Q(0)=Q_{0}$ &Individuals in class $Q$ at $t=0$\\
$H(0)=H_{0}$ & Individuals in class $H$ at $t=0$\\
$\overline{H}(0)=\overline{H}_{0}$ & Individuals in class $\overline{H}$ at time $t=0$\\ 
\bottomrule
\end{tabular}
\end{specialtable}

Model \eqref{eq:modelC} has the following reproduction number:
\begin{equation}
\label{R0}
\mathcal{R}_{0}=\dfrac{\beta a_{2} (1-p)[(l_{H} a_{6} qv + (a_{1} + qv)a_{3})a_{7}
-\delta_{2} \eta_{k} f_{2} (q v + a_{1})]}{a_{0} a_{1} 
\chi (p \phi +a_{2})}=\dfrac{\mathcal{N}}{\mathcal{D}}
\end{equation}
and two equilibrium points, the disease free equilibrium (DFE) point 
\begin{equation}
\label{DFE}
\Sigma_{0}= \left(S_{0},A_{0}, I_{0},Q_{0}, H_{0}, \overline{H}_{0} \right)
=\left(\frac{\Lambda a_{2}}{(p \phi + a_{2})\mu}, 0, 0, 
\frac{p \phi \Lambda}{(p \phi + a_{2}) \mu}, 0, 0 \right)
\end{equation}
and the endemic equilibrium (EE) point that is given by the following:
\begin{equation}
\label{endemic}
\Sigma^{\ast}=\left(S^{\ast}, A^{\ast}, I^{\ast}, Q^{\ast}, 
H^{\ast}, \overline{H}^{\ast} \right)
\end{equation}
with
\begin{equation*}
\begin{split}
S^{\ast}&=\frac{\lambda^{\ast} a_{0} a_{1} a_{2} \chi}{D^{\ast}},         
\quad A^{\ast}= \frac{a_{1} a_{2} \chi \Lambda \lambda^{\ast}}{ D^{\ast}},\\
I^{\ast}&=\frac{\chi \Lambda a_{2} q v \lambda^{\ast} }{D^{\ast}},
\quad Q^{\ast}=\frac{\Lambda((\chi \delta_{1} f_{1}+a_{4} a_{6} a_{7})q v \lambda^{\ast} 
+ a_{0}a_{1} p \phi \chi)}{D^{\ast}},\\
H^{\ast}&=\frac{\Lambda a_{2}a_{6} a_{7} q v \lambda^{\ast}}{ D^{\ast}}, \quad
\overline{H}^{\ast}=\frac{\delta_{2} f_{2} \Lambda a_{2} a_{6} q v \lambda^{\ast}}{D^{\ast}},
\end{split}
\end{equation*}
where $D^{\ast}= ( \chi ( -f_{1} m \omega q v \delta_{1} + a_{0} a_{1} a_{2})-a_{4} a_{6} a_{7} 
m \omega q v) \lambda^{\ast} + \mathcal{D} \mu$. Assuming that the transmission rate is strictly
positive, we have the following:
\begin{equation}
\label{steady}
\lambda^{\ast}=\frac{\beta (A^{\ast} +I^{\ast} +l_{H} H^{\ast})(1-p)}{N^{\ast}}.
\end{equation}

Using \eqref{endemic} in \eqref{steady}, it can be seen 
that the endemic equilibrium satisfies the following:
\begin{equation}
\lambda^{\ast}=\dfrac{(\mathcal{N}-\mathcal{D}) \beta (1-p)}{\mathcal{N} 
+ qv(a_{2} a_{6} (a_{7} (1-l_{H}) + \delta_{2} f_{2}) + \delta_{1} f_{1} \chi
+ a_{4} a_{6} a_{7}) \beta (1-p) }.
\end{equation}

Regarding the stability of the equilibrium points, the following results hold.

\begin{Lemma}[Lemma 3.2 of \cite{PST}]
The disease free equilibrium point $\Sigma_{0}$ \eqref{DFE} 
is locally asymptotically stable if $\mathcal{R}_{0}<1$ and unstable 
if $\mathcal{R}_{0}>1$, where $\mathcal{R}_{0}$ is the basic reproduction 
number \eqref{R0}.
\end{Lemma}

\begin{Lemma}[Lemma 3.3 of  \cite{PST}]
The model \eqref{eq:modelC} has a unique endemic equilibrium 
point whenever $\mathcal{R}_{0} >1.$
\end{Lemma}


\section{Results}
\label{sec3}

In this section, we begin by presenting our proposal for a COVID-19 discrete-time model. 
After that, we show the well-posedness of the model, that is, we prove that the solutions 
of the model are positive and bounded and that the equilibrium points coincide with those 
of the continuous-time model presented in Section~\ref{sec2}. We finalize this section 
by proving the global stability of the disease-free equilibrium point, using a 
suitable Lyapunov function and with the presentation of some numerical simulations, 
using real data, showing the consistency of our theoretical results.


\subsection{The Discrete-Time Model}

One of the important features of the discrete-time epidemic models obtained by the Mickens method is that
they present the same features of the corresponding original continuous-time models. Our 
nonstandard finite discrete difference (NSFD) scheme for solving \eqref{eq:modelC} is 
a numerical dynamically consistent method based on \cite{Mickens05}. Let us define the time 
instants $t_{n} = n h$ with $n$ integer, the step size as $h = t_{n+1} - t_{n}$, 
and $(S_{n}, A_{n}, I_{n},Q_{n}, H_{n}, \overline{H}_{n})$ as the approximated values of  the following:
$$
(S(nh),  A(nh), I(nh),Q(nh), H (nh), \overline{H}(nh)). 
$$

Discretizing system \eqref{eq:modelC} 
using the NSFD scheme, we obtain the following:
\begin{adjustwidth}{-4.6cm}{0cm}
\begin{equation}
\label{eq:modelD}
\begin{cases}
\dfrac{S_{n+1}-S_{n}}{\psi(h)} = \Lambda + \omega m Q_{n+1} 
- \left[ \lambda_{n} (1-p) + \phi p + \mu\right] S_{n+1},\\[0.2 cm]
\dfrac{A_{n+1}-A_{n}}{\psi(h)}= \lambda_{n}(1-p) S_{n+1} - (q v +\mu) A_{n+1}, \\[0.2 cm]
\dfrac{I_{n+1}-I_{n}}{\psi(h)} = q v A_{n+1} - ( \delta_{1} + \mu)I_{n+1},\\[0.2 cm]
\dfrac{Q_{n+1}-Q_{n}}{\psi(h)} = \phi p S_{n+1} + \delta_{1} f_{1} I_{n+1} 
+ \delta_{2}(1-f_{2}-f_{3})H_{n+1} - (\omega m + \mu) Q_{n+1},\\[0.2 cm]
\dfrac{H_{n+1}-H_{n}}{\psi(h)}=\delta_{1}(1-f_{1})I_{n+1}+\eta (1-\kappa)\overline{H}_{n+1}
-\left(\delta_{2}(1-f_{2}-f_{3})+\delta_{2}f_{2}+ \alpha_{1}f_{3}+\mu\right)H_{n+1}, \\[0.2 cm]
\dfrac{\overline{H}_{n+1}-\overline{H}_{n}}{\psi(h)}=\delta_{2}f_{2} H_{n+1} 
- \left(\eta(1-\kappa) + \alpha_{2} \kappa + \mu\right)\overline{H}_{n+1},
\end{cases}
\end{equation} 
\end{adjustwidth}
where the denominator function is $\psi(h)=\frac{exp(\mu h)-1}{\mu}$ \cite{Mickens07}. 
Throughout this work, for brevity, we write $\psi(h)=\psi$.

\begin{Remark}
In the continuous model, the reproduction number \eqref{R0} can be rewritten as the following:
\begin{equation}
R_{0}=\dfrac{\beta a_{2} (1-p) [l_{H} a_{6} a_{7} q v + (a_{1} + q v) \chi]}{ a_{0} a_{1} 
\chi (p \phi +a_{2})}=\frac{\beta a_{2} (1-p) [l_{H} a_{6} a_{7} q v + (a_{0} 
+\delta_{1}) \chi]}{a_{0} a_{1} \chi (p \phi +a_{2})}.
\end{equation}
\end{Remark}

Let us consider the region $\Omega=\left\{(S,A,I,Q,H,\overline{H}) 
\in (\mathbb{R}_{0}^{+})^{6} : 0 < N \le \frac{\Lambda}{\mu}\right\}$. 
Our next lemma shows that the feasible region is equal to the continuous case.

\begin{Lemma}
Any solution of $(S_{n},A_{n},I_{n},Q_{n}, H_{n},\overline{H}_{n})$ 
of model \eqref{eq:modelD} with positive initial conditions 
is positive and ultimately bounded in $\Omega$.
\end{Lemma}

\begin{proof}
Since model \eqref{eq:modelD} is linear in 
$(S_{n},A_{n},I_{n},Q_{n}, H_{n},\overline{H}_{n})$, we can rewrite it as the following:
\begin{equation}
\label{eq:modelED}
\begin{cases}
S_{n+1} = \dfrac{ \Lambda \psi + \omega m \psi Q_{n+1} 
+ S_{n}}{1+ [ \lambda_{n} (1-p) + \phi p + \mu] \psi},\\[0.4 cm]
A_{n+1}=\dfrac{A_{n}+ \lambda_{n}(1-p)\psi S_{n+1}}{1 + (q v +\mu)\psi }, \\[0.4 cm]
I_{n+1}=\dfrac{I_{n}+ q v \psi A_{n+1}}{1+  ( \delta_{1} + \mu)\psi},\\[0.4 cm]
Q_{n+1} = \dfrac{Q_{n} + \psi( \phi p S_{n+1} + \delta_{1} f_{1} I_{n+1} 
+ \delta_{2}(1-f_{2}-f_{3})H_{n+1})}{1+  (\omega m + \mu)\psi},\\[0.4 cm]
H_{n+1}=\dfrac{H_{n}+ \psi(\delta_{1}(1-f_{1})I_{n+1}+\eta (1-\kappa)\overline{H}_{n+1})}{1
+(\delta_{2}(1-f_{2}-f_{3})+\delta_{2}f_{2}+ \alpha_{1}f_{3}+\mu)\psi}\\[0.4 cm]
\overline{H}_{n+1}=\dfrac{\overline{H}_{n}+\delta_{2}f_{2} 
\psi H_{n+1}}{1 + (\eta(1-\kappa) + \alpha_{2} \kappa + \mu)\psi}.
\end{cases}
\end{equation}
 
Since all parameters of model \eqref{eq:modelED} are positive and the initial conditions 
are also positive, then, by induction, $S_{n}\ge 0$, $A_{n}\ge 0$, $I_{n}\ge 0$, $Q_{n}\ge 0$, $H_{n}\ge 0$ 
and $\overline{H}_{n}\ge 0$ for all $n \in \mathbb{N}$. Regarding the boundedness of the solutions, 
the total population $N_{n}=S_{n}+A_{n}+I_{n}+Q_{n}+H_{n}+\overline{H}_{n}$. 
Adding all equations of \eqref{eq:modelED}, we obtain the following:
\begin{align*}
\dfrac{N_{n+1}-N_{n}}{\psi}&=\Lambda -\mu N_{n+1}- \alpha_{1} 
f_{3} H_{n+1} - \alpha_{2} \kappa \overline{H}_{n+1}\le \Lambda - \mu N_{n+1}\\
\Leftrightarrow N_{n+1} &\le \dfrac{\Lambda \psi}{1+\mu \psi}+ \dfrac{N_{n}}{1+\mu \psi}.
\end{align*}

By Lemma~2.2 of Shi and Dong \cite{S},
\begin{equation*}
N_{n} \le \dfrac{\Lambda}{\mu} + \left( \dfrac{1}{1 + \mu \psi}\right)^{n}\left( 
\dfrac{\Lambda}{\mu} - N_{0}\right)
\end{equation*}
so, if $N_{0} \le \frac{\Lambda}{\mu}$, then $N_{n} \le \frac{\Lambda}{\mu}$ for all $n \in \mathbb{N}$. 
Thus, $\Omega$ is the biologically feasible region.
\end{proof}

Solving $X^{\ast}=F(X^{\ast})$ in system \eqref{eq:modelED} 
we can see that there exists two equilibrium points:
\begin{itemize}
\item The disease free equilibrium (DFE) point
\begin{equation}
\label{E0}
E_{0}=\left( \dfrac{\Lambda a_{2}}{ \mu (p \phi +a_{2})},0, 
\dfrac{\Lambda p \phi}{\mu (p \phi +a_{2})}, 0, 0,0 \right)
=(S_{0}, A_{0}, Q_{0}, I_{0}, H_{0}, \overline{H}_{0});
\end{equation}
\item The endemic equilibrium (EE) point
\begin{equation}
E^{\ast}=(S^{\ast}, A^{\ast}, I^{\ast}, Q^{\ast}, H^{\ast}, \overline{H}^{\ast}),
\end{equation}
where
\begin{equation}
\label{EE}
\begin{gathered}
S^{\ast} =\dfrac{a_{0} a_{1} a_{2} \Lambda \chi}{D_{A}},
\quad A^{\ast}=\dfrac{a_{1} a_{2} \Lambda \chi \lambda^{\ast}(1-p)}{D_{A}}, \\
I^{\ast}=\dfrac{a_{2} q v \Lambda \chi \lambda^{\ast}(1-p)}{D_{A}},
\quad Q^{\ast} = \dfrac{\Lambda (a_{0} a_{1} p \phi \chi + \lambda^{\ast} (1-p)  
q v (\chi f_{1} \delta_{1} + a_{4} a_{6} a_{7}))}{D_{A}},\\
H^{\ast}=\dfrac{a_{2} a_{6} a_{7} q v \Lambda \lambda^{\ast}(1-p)}{D_{A}}, 
\quad \overline{H}^{\ast}=\dfrac{a_{2} a_{6} \delta_{2} f_{2} q v  \Lambda \lambda^{\ast} (1-p)}{D_{A}},
\end{gathered}
\end{equation} 
\end{itemize}
and $D_{A}= \mu \mathcal{D} + \lambda^{\ast} (1-p) ( \chi ( a_{2} a_{1} a_{0} 
-f_{1} \delta_{1} \omega m q v) - a_{4} a_{6} a_{7} q v w m)$. In this point, 
\begin{equation}
\lambda^{\ast}=\dfrac{\mathcal{D} (\mathcal{R}_{0} -1)}{(1-p)(\mathcal{N} 
+ q v ( a_{2} a_{6} ( a_{7} (1-l_{H}) + \delta_{2}f_{2}) 
+ \chi f_{1} \delta_{1} + a_{4} a_{6} a_{7}))}.
\end{equation}


\subsection{Global Stability}

We now prove the global stability of \eqref{E0}.

\begin{Theorem}
For the discretized system, the DFE point $E_{0}$ is globally stable if 
$\mathcal{R}_{0} <1$. If $\mathcal{R}_{0} >1$, then the DFE is unstable.
\end{Theorem}

\begin{proof}
Let us define the discrete Lyapunov function $L_{n}$ as the following:
\begin{equation}
L_{n}(S_{n}, A_{n}, I_{n},Q_{n},H_{n}, \overline{H}_{n})
=\frac{1}{\psi}\left[ S_{0} G\left( \frac{S_{n}}{S_{0}} \right) 
+ A_{n} + I_{n} + Q_{0} G\left( \frac{Q_{n}}{Q_{0}}\right) + H_{n} 
+ \overline{H}_{n} \right].
\end{equation}

Hence, $L_{n}(S_{n}, A_{n}, I_{n},Q_{n},H_{n}, \overline{H}_{n}) \ge 0$ 
for all $\xi_{n}  \ge 0$, $\xi \in \left\{S, A, I,Q,H, \overline{H} \right\}$. 
Moreover, $L_{n}(S_{n}, A_{n}, I_{n},Q_{n},H_{n}, \overline{H}_{n}) =0$ 
if and only if $(S_{n}, A_{n}, I_{n},Q_{n},H_{n}, \overline{H}_{n})=E_{0}$.
Computing $\Delta L_{n}=L_{n+1}-L_{n}$, we have the following:
\begin{align*}
&\Delta L_{n}=\frac{1}{\psi}\left[ S_{0} \left(G\left( \frac{S_{n+1}}{S_{0}}\right)
-G\left( \frac{S_{n}}{S_{0}}\right) \right) +(A_{n+1}-A_{n}) +(I_{n+1}-I_{n})\right]\\
& +  \frac{1}{\psi}\left[Q_{0} \left(G\left( \frac{Q_{n+1}}{Q_{0}}\right)
-G\left( \frac{Q_{n}}{Q_{0}}\right)\right) +(H_{n+1}-H_{n})+ (\overline{H}_{n+1}
-\overline{H}_{n})\right],
\end{align*}
where $G(x)= x - \ln(x) -1$. Note that $G(x) \ge 0$, for all $x \ge 0$, and $G(x)=0$ 
if and only if $x=1$. For each $\xi \in \left\{S, A, I,Q,H, \overline{H} \right\}$,
\begin{equation*}
G\left(\frac{\xi_{n+1}}{\xi_{0}}\right)-G\left(\frac{\xi_{n}}{\xi_{0}}\right) 
\le \frac{1}{\xi_{0}} \left( 1- \frac{\xi_{0}}{\xi_{n+1}} \right)\left( \xi_{n+1}-\xi_{n}\right).
\end{equation*}

Therefore,
\begin{align*}
&\Delta L_{n} \le \frac{1}{\psi}\left[ \left( 1 - \frac{S_{0}}{S_{n+1}}\right)(S_{n+1}-S_{n}) 
+ (A_{n+1}-A_{n})+(I_{n+1}+I_{n})\right]\\
&\quad +\frac{1}{\psi}\left[\left(1-\frac{Q_{0}}{Q_{n+1}}\right)(Q_{n+1}-Q_{n})
+(H_{n+1}-H_{n})+(\overline{H}_{n+1}-\overline{H}_{n})\right].
\end{align*}

From equations of system \eqref{eq:modelD}, we obtain the following:
\begin{align*}
&\Delta L_{n}\le \left( 1-\frac{S_{0}}{S_{n+1}}\right)(\Lambda + \omega m Q_{n+1}
- (\lambda_{n} (1-p) + \phi p + \mu)S_{n+1} )\\
&+(\lambda_{n} (1-p)S_{n+1}-(q v +\mu)A_{n+1}) + (q v A_{n+1} - (\delta_{1}+ \mu)I_{n+1})\\
&+\left( 1-\frac{Q_{0}}{Q_{n+1}}\right)(\phi p S_{n+1}+\delta_{1}f_{1}I_{n+1}
+\delta_{2}(1-f_{2}-f_{3})H_{n+1} - (\omega m + \mu) Q_{n+1} ) \\
&+(\delta_{1}(1-f_{1}) I_{n+1} + \eta (1-\kappa) \overline{H}_{n+1}
-\left[ \delta_{2}(1-f_{2}-f_{3})+\delta_{2}f_{2} +\alpha_{1} f_{3} + \mu\right] H_{n+1})\\
&+ \delta_{2} f_{2} H_{n+1} - (\eta(1-\kappa) + \alpha_{2} \kappa + \mu) \overline{H}_{n+1}
\end{align*}
and, at the disease-free equilibrium point $E_{0}$, we have the following relations:
\begin{equation}
\begin{cases}
0= \Lambda + \omega m Q_{0} - (\phi p + \mu) S_{0}, \\
0=\phi p S_{0} - (\omega m + \mu) Q_{0}, 
\end{cases} 
\Leftrightarrow
\begin{cases}
\Lambda=  - \omega m Q_{0} + (\phi p + \mu) S_{0} ,\\
\phi p S_{0}= (\omega m + \mu) Q_{0}. 
\end{cases}
\end{equation}

Substituting in $\Delta L_{n}$, and simplifying, we obtain the following:
\begin{align*}
&\Delta L_{n} \le -\frac{(\phi p + \mu)}{S_{n+1}} (S_{n+1}-S_{0})^{2} 
- \mu (A_{n+1}+ I_{n+1}) - \left( \alpha_{1} f_{3} + \mu\right) H_{n+1} \\
&-\frac{Q_{0}}{Q_{n+1}} \left(  \delta_{1} f_{1} I_{n+1} 
+ \delta_{2} (1-f_{2} -f_{3})H_{n+1}\right)-(\alpha_{2} \kappa + \mu) \overline{H}_{n+1}
+\lambda_{n} (1-p)S_{0}\\
&+ \omega m Q_{0}  \left(  1-\frac{S_{0}}{S_{n+1}}\right) \left( \frac{Q_{n+1}}{Q_{0}}-1\right)
+\phi p S_{0} \left( 1-\frac{Q_{0}}{Q_{n+1}}\right)\left(
\frac{S_{n+1}}{S_{0}}-\frac{Q_{n+1}}{Q_{0}}\right).
\end{align*}

Let 
\begin{equation*}
F(I,Q,H, \overline{H}) = \left( \alpha_{1} f_{3} + \mu\right) H
+\frac{Q_{0}}{Q} \left(  \delta_{1} f_{1} I + \delta_{2} (1-f_{2} -f_{3})H\right)
+(\alpha_{2} \kappa + \mu) \cdot \overline{H}.
\end{equation*}

Then,
\begin{align}
\label{ineqL}
&\Delta L_{n} \le -\frac{(\phi p + \mu)}{S_{n+1}} (S_{n+1}-S_{0})^{2} 
- \mu (A_{n+1}+ I_{n+1})-F(I_{n+1},Q_{n+1},H_{n+1}, \overline{H}_{n+1}) \nonumber \\
&+ \omega m Q_{0}  \left(  1-\frac{S_{0}}{S_{n+1}}\right) \left( \frac{Q_{n+1}}{Q_{0}}-1\right)
+\phi p S_{0} \left( 1-\frac{Q_{0}}{Q_{n+1}}\right)\left(
\frac{S_{n+1}}{S_{0}}-\frac{Q_{n+1}}{Q_{0}}\right)\\
&+\lambda_{n} (1-p)S_{0}.\nonumber
\end{align}

From the equations of \eqref{eq:modelED}, it can be seen that
\begin{align}
\label{ineqlambda}
\lambda_{n} (1-p) S_{0} & \le \frac{\beta \Lambda (1-p) a_{2}}{ \mu (p \phi 
+ a_{2})}(A_{n}+ I_{n}+ l_{H} H_{n}) \nonumber\\
&  \le \frac{\beta \Lambda (1-p) a_{2}}{ \mu (p \phi + a_{2})}\left(  A_{n}
+ \frac{q v}{a_{1}} A_{n} + \frac{a_{6} a_{7} q v}{ \chi a_{1}} l_{H} A_{n}\right)\nonumber\\
&\le \frac{\Lambda \mathcal{N}}{ \mu (p \phi + a_{2})} A_{n}
\end{align}
and 
\begin{align}
\label{ineqDFE}
&\omega m Q_{0}  \left(  1-\frac{S_{0}}{S_{n+1}}\right) \left( \frac{Q_{n+1}}{Q_{0}}-1\right)
+\phi p S_{0} \left( 1-\frac{Q_{0}}{Q_{n+1}}\right)\left(
\frac{S_{n+1}}{S_{0}}-\frac{Q_{n+1}}{Q_{0}}\right) \nonumber \\
&\le  w m Q_{0}\left(  g \left( \frac{Q_{n+1}}{Q_{0}}\right) 
+ g\left( \frac{S_{0}}{S_{n+1}} \right)- g\left( \frac{S_{0} Q_{n+1}}{S_{n+1} Q_{0}}\right) \right)\\
&\quad +\phi p S_{0} \left(  -g\left( \frac{Q_{n+1}}{Q_{0}}\right) - g\left( \frac{Q_{0} S_{n+1}}{Q_{n+1}
S_{0}}\right)+g\left( \frac{S_{n+1}}{S_{0}}\right)   \right).\nonumber
\end{align}

Since
\begin{equation*}
w m Q_{0} < (w m + \mu) Q_{0}=a_{2}Q_{0}=\phi p S_{0},
\end{equation*}
\eqref{ineqDFE} is equal to the following:
\begin{align}
\label{ineqDFE1}
&g \left( \frac{Q_{n+1}}{Q_{0}}\right) \left( w m Q_{0}-\phi p S_{0}\right)+w m Q_{0} g\left(
\frac{S_{0}}{S_{n+1}} \right)+\phi p S_{0} g\left( \frac{S_{n+1}}{S_{0}}\right)\nonumber\\
&-w m Q_{0} g\left( \frac{S_{0} Q_{n+1}}{S_{n+1} Q_{0}}\right) -\phi p S_{0} 
g\left( \frac{Q_{0} S_{n+1}}{Q_{n+1} S_{0}}\right)  \\
\le &-w m Q_{0} g\left( \frac{S_{0} Q_{n+1}}{S_{n+1} Q_{0}}\right) -\phi p S_{0} g\left( 
\frac{Q_{0} S_{n+1}}{Q_{n+1} S_{0}}\right) 
+\phi p S_{0}\frac{(S_{n+1}-S_{0})^{2}}{S_{0}S_{n+1}} \nonumber\\
&+g \left( \frac{Q_{n+1}}{Q_{0}}\right) \left( w m Q_{0}-\phi p S_{0}\right). \nonumber 
\end{align}

Gathering the information from  \eqref{ineqlambda} and \eqref{ineqDFE1},  
\eqref{ineqL} becomes the following:
\begin{align*}
\Delta L_{n} &\le -\frac{\mu}{S_{n+1}} (S_{n+1}-S_{0})^{2} -w m Q_{0} g\left( 
\frac{S_{0} Q_{n+1}}{S_{n+1} Q_{0}}\right) -\phi p S_{0} 
g\left( \frac{Q_{0} S_{n+1}}{Q_{n+1} S_{0}}\right)\\ 
&-g \left( \frac{Q_{n+1}}{Q_{0}}\right) \left( -w m Q_{0}+\phi p S_{0}\right)
-F(I_{n+1},Q_{n+1},H_{n+1}, \overline{H}_{n+1}) \\
&+\frac{\Lambda \mathcal{N}}{ \mu (p \phi + a_{2})} A_{n} - \mu (A_{n+1}+ I_{n+1}).
\end{align*}

Once more, from system \eqref{eq:modelED} it can seen that 
$- \mu (A_{n+1}+ I_{n+1}) < -a_{0}A_{n}$, so
\begin{align*}
&\frac{\Lambda \mathcal{N}}{ \mu (p \phi + a_{2})} A_{n} - \mu (A_{n+1}+ I_{n+1})
< \frac{A_{n} \mathcal{D}}{a_{1}\chi (p \phi + a_{2})} (\mathcal{R}_{0}-1).
\end{align*}

Therefore,
\begin{align*}
\Delta L_{n}& \le  -\frac{\mu}{S_{n+1}} (S_{n+1}-S_{0})^{2} -w m Q_{0} g\left( \frac{S_{0}
Q_{n+1}}{S_{n+1} Q_{0}}\right) -\phi p S_{0} g\left( \frac{Q_{0} S_{n+1}}{Q_{n+1} S_{0}}\right)\\ 
&-g \left( \frac{Q_{n+1}}{Q_{0}}\right) \left( -w m Q_{0}+\phi p S_{0}\right)-F(I_{n+1},Q_{n+1},
H_{n+1}, \overline{H}_{n+1})\\
&+\frac{A_{n} \mathcal{D}}{a_{1}\chi (p \phi + a_{2})} (\mathcal{R}_{0}-1).
\end{align*}

Hence, if $\mathcal{R}_{0} <1$, then $\Delta L_{n}\le 0$ for all $n \ge 0$, 
that is, $L_{n}$ is a monotone decreasing sequence. We have $L_{n} \ge 0$. 
Then there is a limit for $\underset{n \to \infty}{\lim} L_{n} \ge 0$. 
Therefore, $\underset{n \to \infty}{\lim} \Delta L_{n}=0$ implies that 
$\underset{n \to \infty}{\lim}S_{n}=S_{0}$, $\underset{n \to \infty}{\lim}Q_{n}=Q_{0}$, 
$\underset{n \to \infty}{\lim}A_{n}= \underset{n \to \infty}{\lim}I_{n}
=\underset{n \to \infty}{\lim}H_{n}=\underset{n \to \infty}{\lim}\overline{H}_{n}=0$.
So, if $\mathcal{R}_{0} < 1$, then $E_{0}$ is globally  asymptotically stable. 
\end{proof}


\subsection{Numerical Simulations}

In this section, we show, numerically, that our discretized model describes well 
the transmission dynamics of COVID-19 in Portugal, from 2 March to 4  May 2020. 
Our data were obtained from the daily reports from DGS, available in \cite{DGS}, 
and since we consider the spread from 2 March, the initial time $t=0$ 
corresponds to 2 March 2020.

The initial values are the same as those used in \cite{PST}. Regarding the parameters, 
using \cite{Pordata}, we can obtain the values from 2019, so the parameters 
were updated with respect to \cite{PST}, but the difference is very small. 
All the values used in our simulations are presented in Table~\ref{tab2}.

\begin{specialtable}[H] 
\caption{Parameter values and initial conditions of \eqref{eq:modelC}.\label{tab2}}
\begin{tabular}{lm{9.05cm}<{\raggedright}l}
\toprule
\textbf{Parameter}	& \textbf{Value} & \textbf{Reference} \\
\midrule
$\Lambda$ & (86,579 + 26,080)/365(person day$^{-1}$)	& \cite{Pordata} \\
$\mu$     & 111,793/(365 $\times N_{0}$)(day$^{-1}$)& \cite{Pordata}\\
$\beta$   & 1.93 (day$^{-1}$)& \cite{PST} \\
$l_{A}$   & 1 (dimensionless)& \cite{PST}\\
$l_{H}$   & 0.1 (dimensionless)& \cite{PST}\\
$\phi$    & 1/12& \cite{RP}\\
$\nu$     & 1/5& \cite{Who}\\
$\delta_{1}$ & 1/3 (day$^{-1}$)& \cite{PST}\\
$\delta_{2}$ & 1/3 (day$^{-1}$)& \cite{PST}\\
$\eta$       & 1/7 (day$^{-1}$)& \cite{PST}\\
$\omega$     & 1/31 (day$^{-1}$)& \cite{PST}\\
$\alpha_{1}$ & 1/7 (day$^{-1}$)& \cite{PST}\\
$\alpha_{2}$ & 1/15 (day$^{-1}$)& \cite{PST}\\
$p$ & 0.674  & \cite{Pordata,Negocios}\\
$q$ & 0.15   & \cite{Noticias}\\
$f_{1}$ & 0.96  & \cite{DGS}\\
$f_{2}$ & 0.21  & \cite{DGS}\\
$f_{3}$ & 0.03  & \cite{DGS}\\
$\kappa$ & 0.03 & \cite{PST}\\
$m$ & 0.075& \cite{PST}\\
$S_{0}$& 10,286,285 (person)& \cite{DGS,Pordata,Who,Noticias}\\
$A_{0}$ & 13 (person)& \cite{DGS,Who,Noticias}\\
$I_{0}$ & 2 (person)&\cite{DGS}\\
$Q_{0}$ & 0 (person)& \cite{PST}\\
$H_{0}$ & 0 (person)& \cite{DGS} \\
$\overline{H}_{0}$ & 0 (person) & \cite{DGS}\\ 
$D_{0}$ & 0 (person)& \cite{DGS}\\
\bottomrule
\end{tabular}
\end{specialtable}

\begin{figure}[H]
\includegraphics[width=10.5 cm]{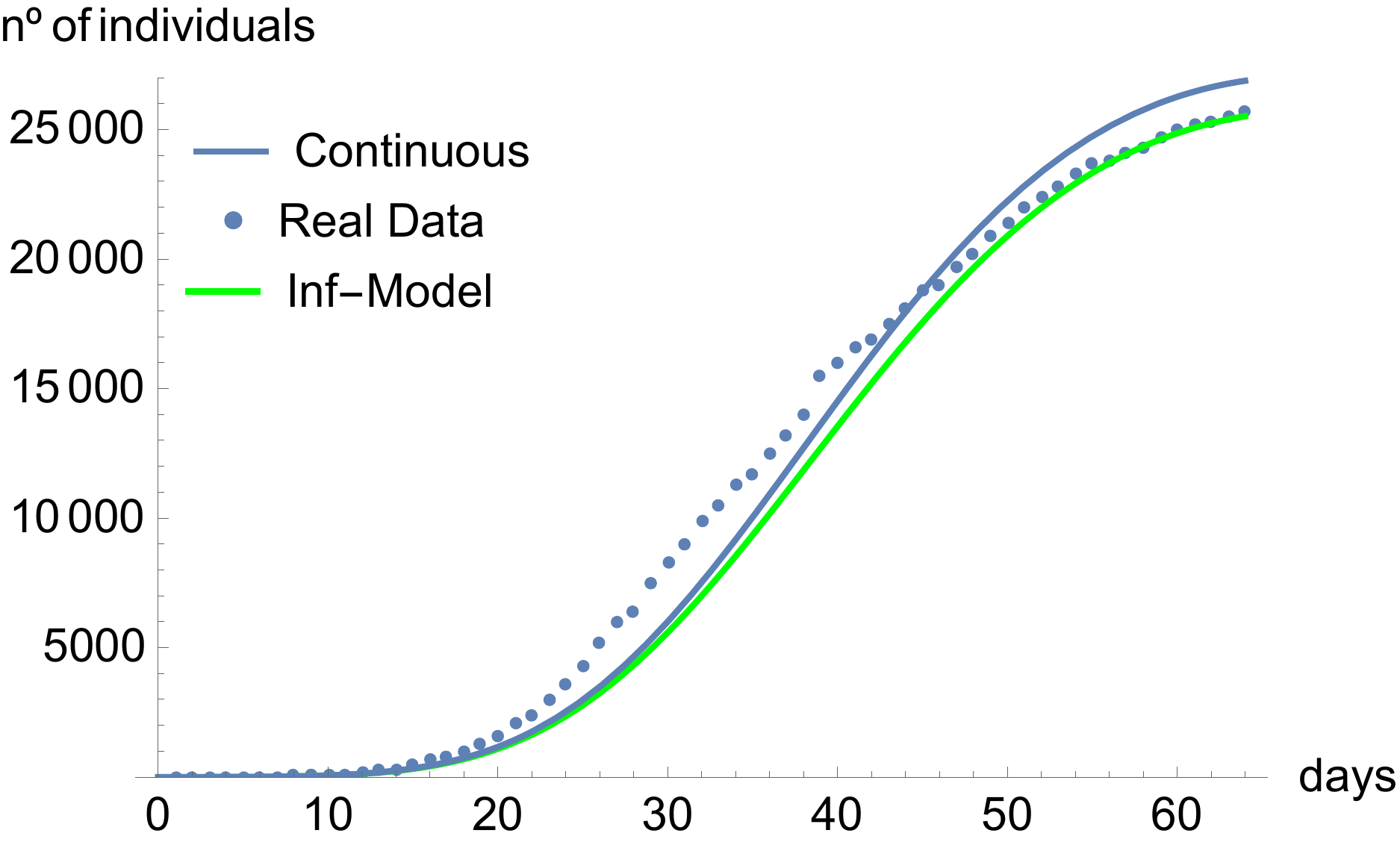}
\caption{Number of infected individuals predicted by the discrete-time model: 
solid green line. Real data: dotted line. Number of infected individuals 
predicted by the continuous model: solid blue line. \label{fig1}}
\end{figure}

Using the values of Table~\ref{tab2}, in Figure~\ref{fig1}, we compare the number 
of infected individuals predicted by our discrete-time model, the ones predicted 
by the continuous model of \cite{PST}, and the real data.

\begin{figure}[H]
\includegraphics[width=9 cm]{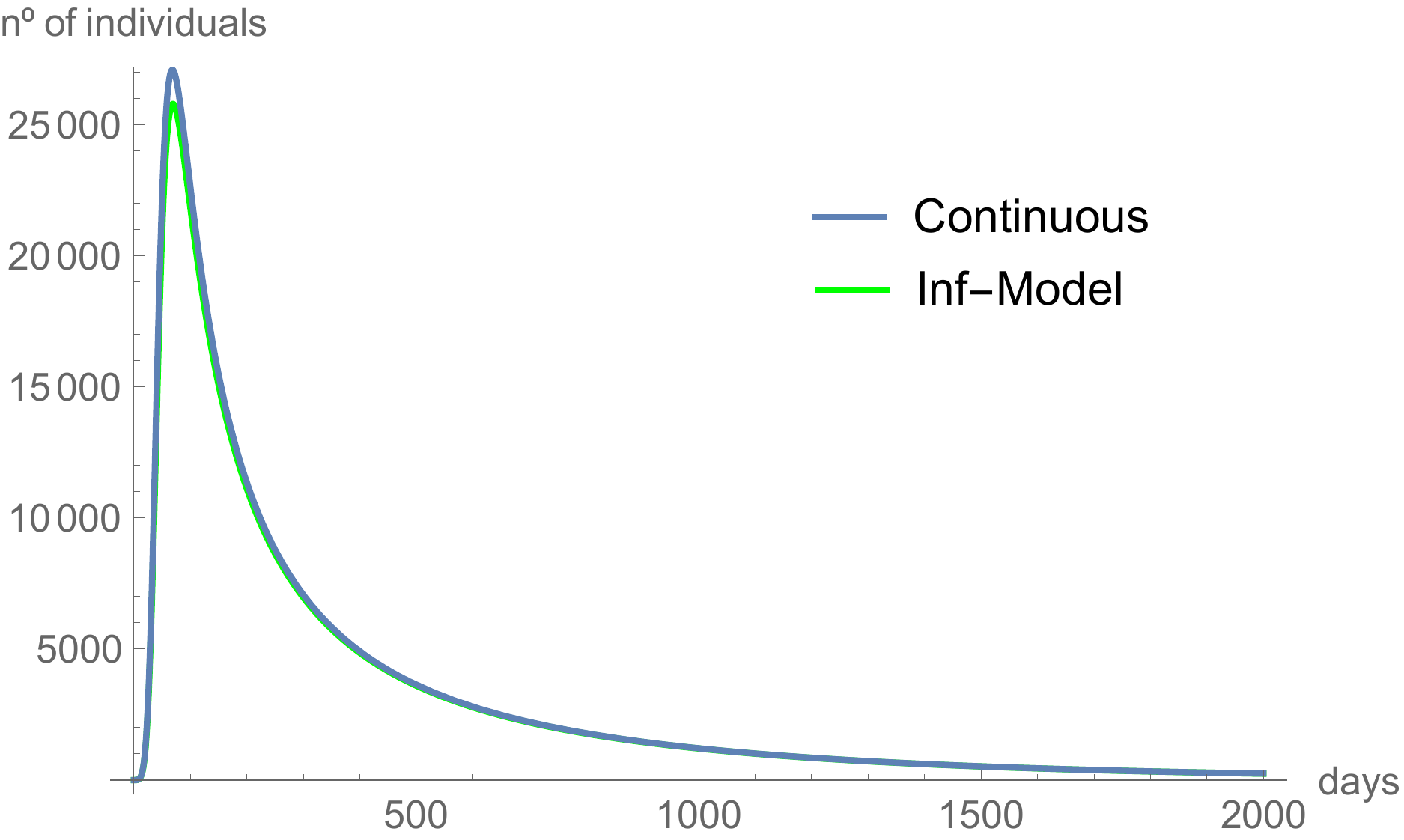}
\caption{The infected individuals tend to disappear with time.\label{fig2}}
\end{figure}

Figure~\ref{fig2} shows the consistency of our result. We compare the predictions 
from the continuous model  with the predictions from the discrete-time model. 
Using the parameter values of Table~\ref{tab2}, $R_{0} \approx 0.95$, 
the number of infected and asymptomatic individuals, as well as the ones in the hospital 
and intensive care units, tend to disappear. The susceptible individuals 
and those in quarantine tend to the equilibrium point $E_{0}$ as time increases.

\begin{figure}[H]
\includegraphics[width=9 cm]{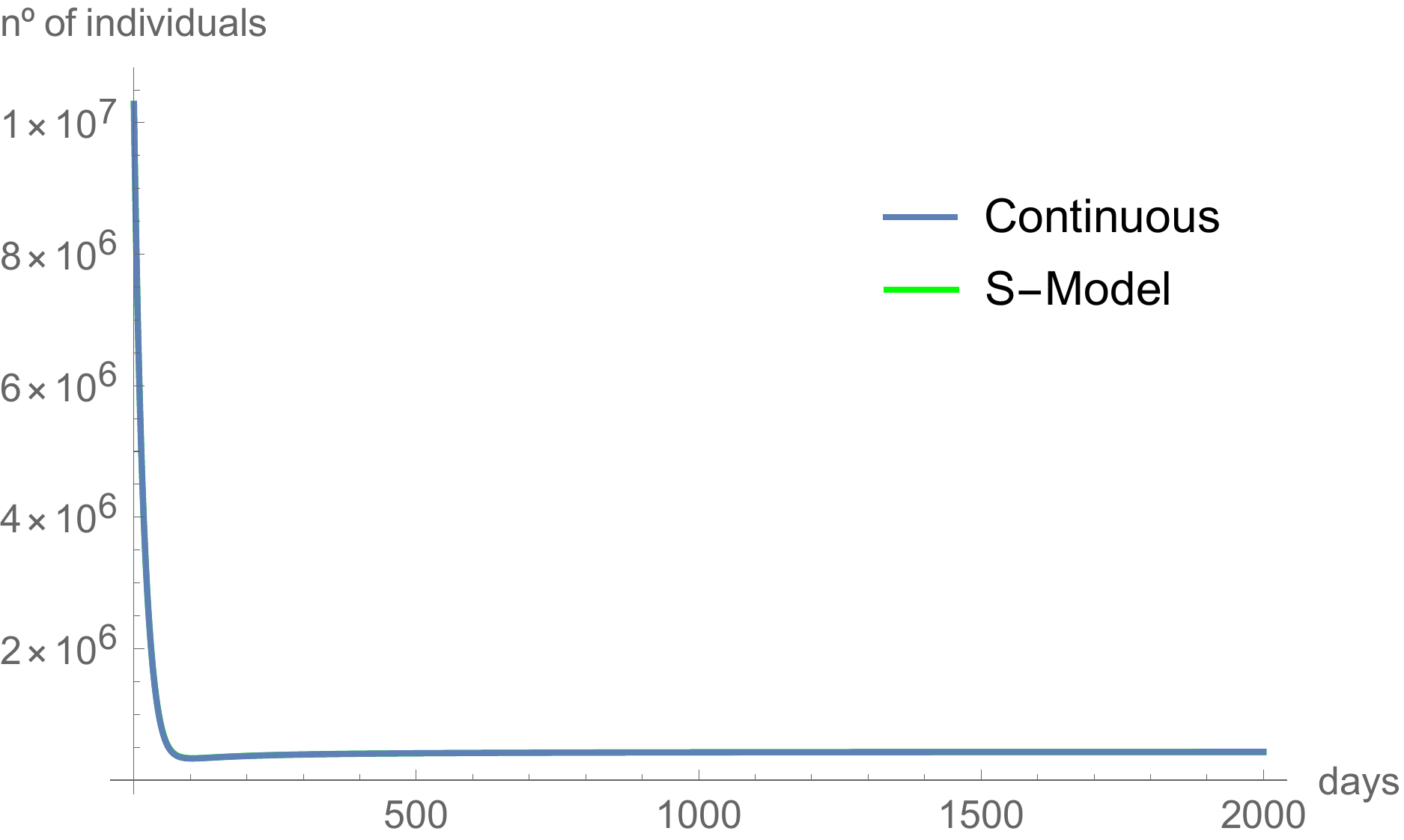}
\caption{The susceptible individuals tend quickly to $S_{0}$. \label{fig3}}
\end{figure}

\begin{figure}[H]
\includegraphics[width=9 cm]{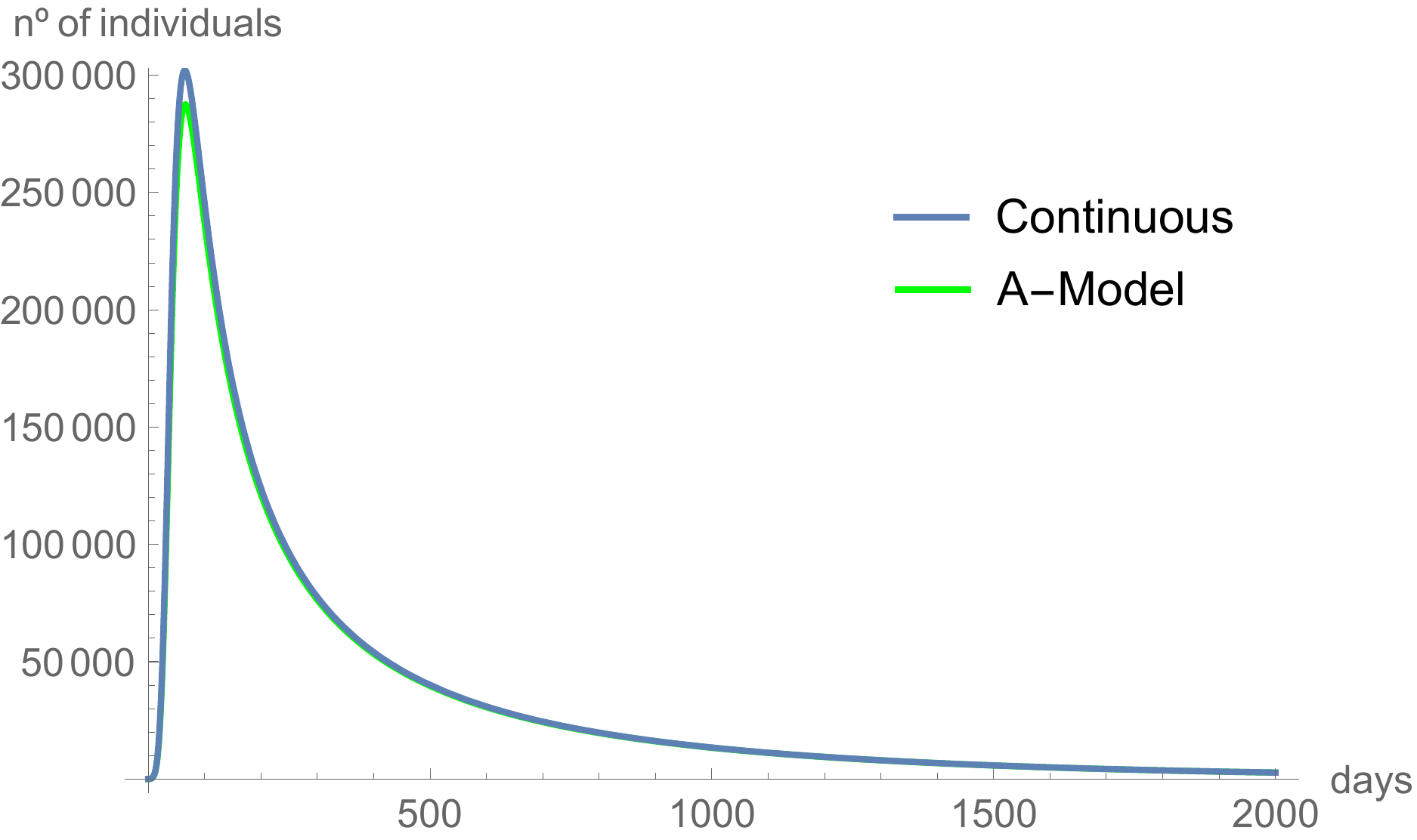}
\caption{The asymptomatic individuals tend to disappear with time. \label{fig4}}
\end{figure}

\begin{figure}[H]
\includegraphics[width=9 cm]{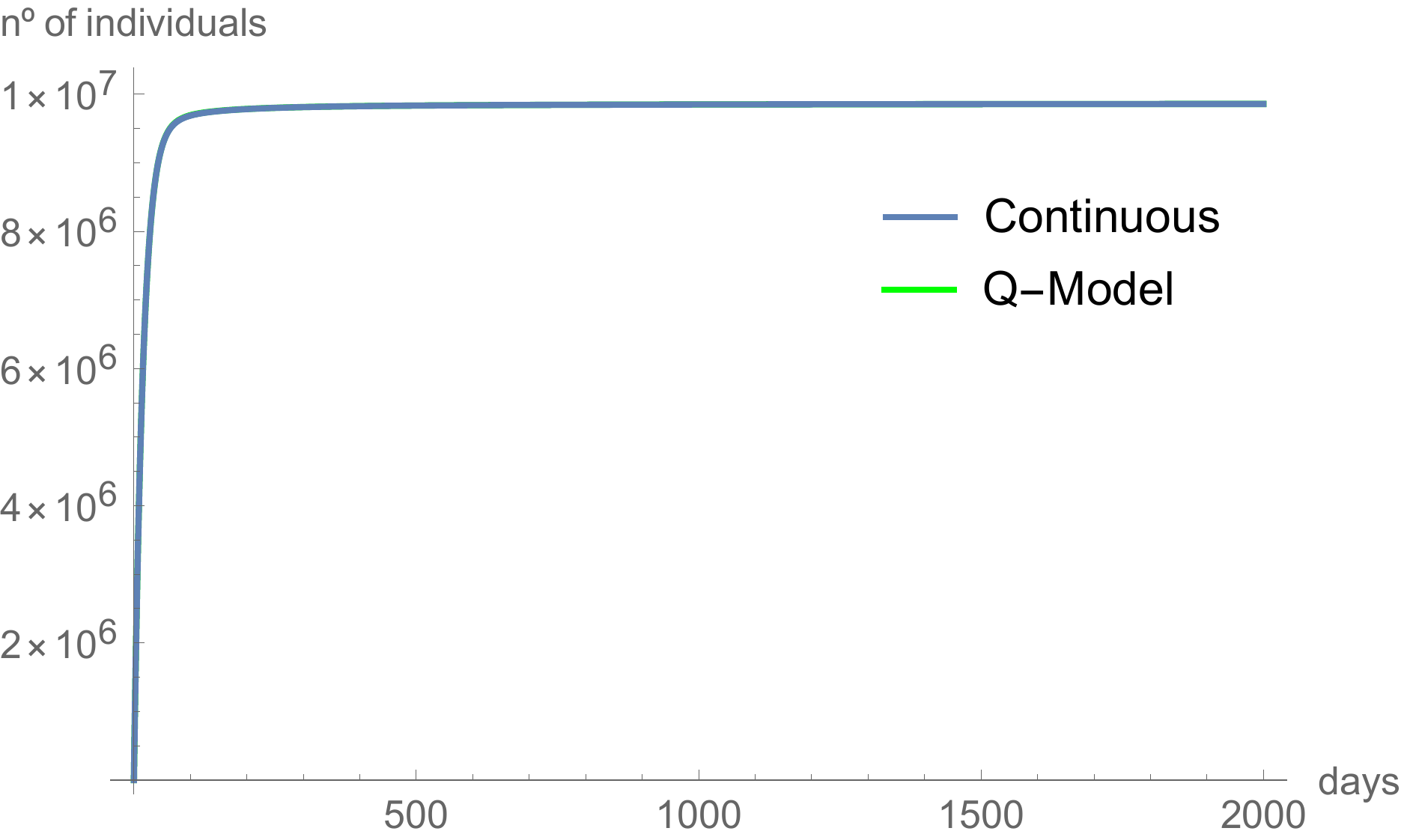}
\caption{The number of individuals  in  quarantine tend quickly to $Q_{0}$. \label{fig5}}
\end{figure}

\begin{figure}[H]
\includegraphics[width=9 cm]{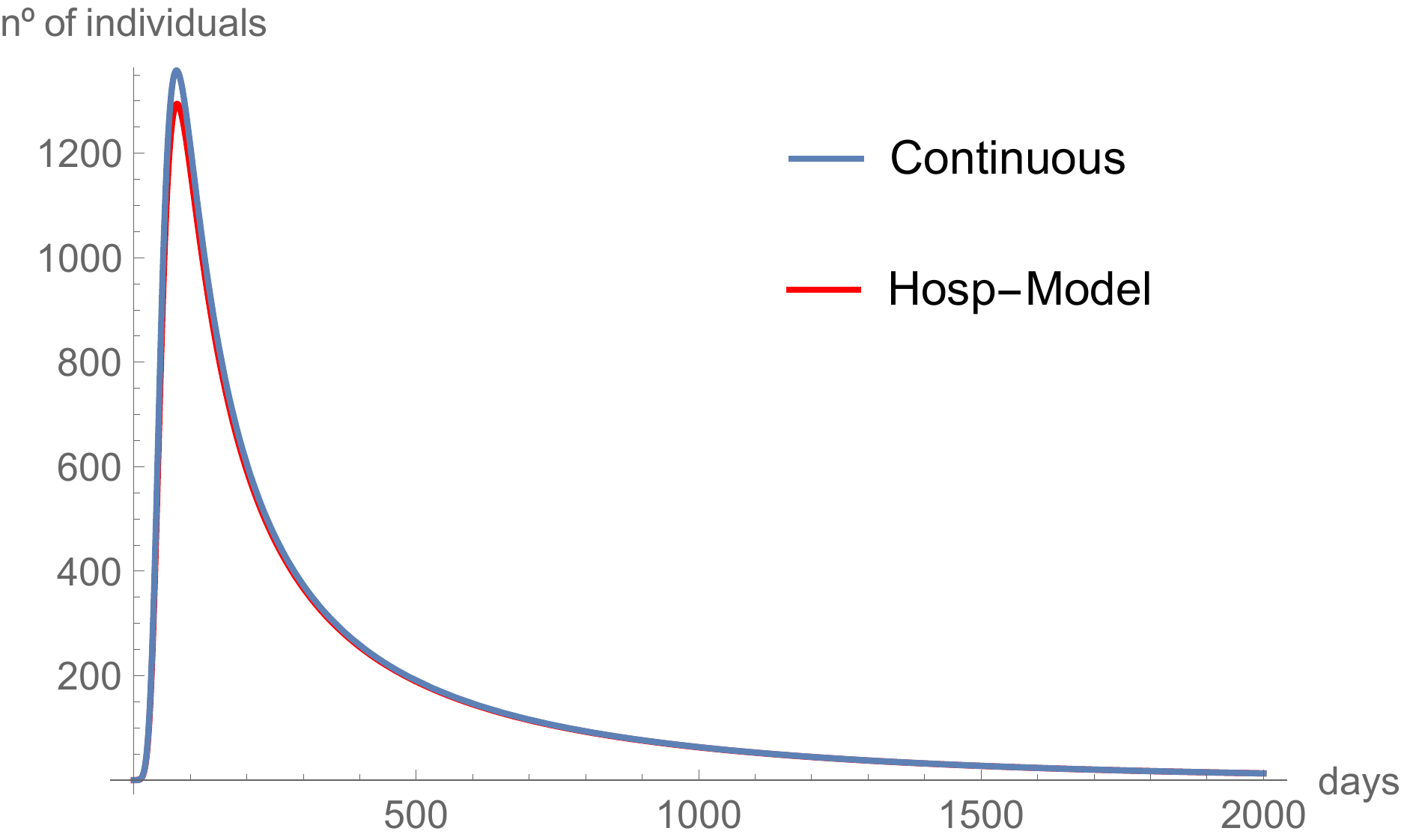}
\caption{The hospitalized individuals tend to disappear with time. \label{fig6}}
\end{figure}

\begin{figure}[H]
\includegraphics[width=9 cm]{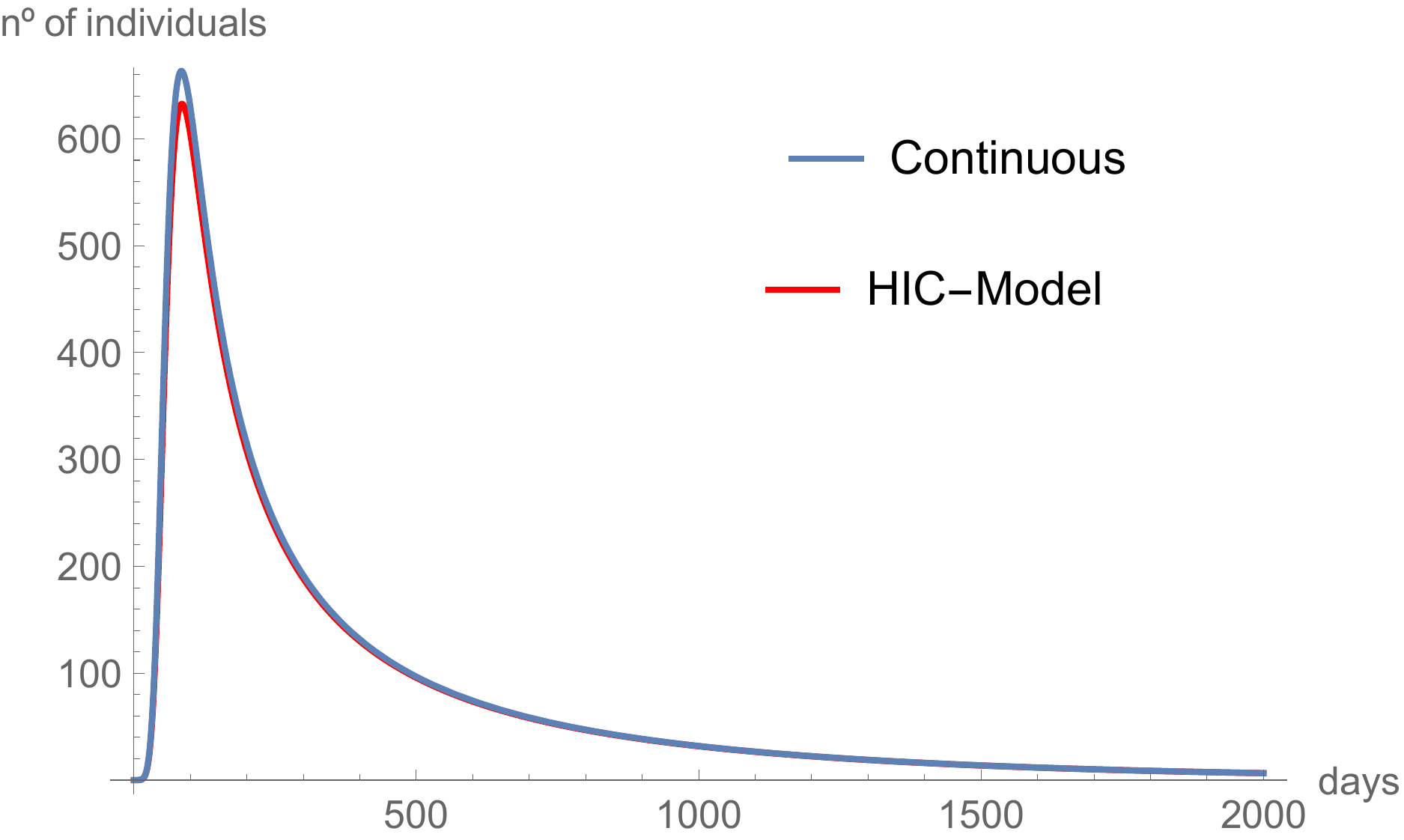}
\caption{The individuals  in intensive care also tend to disappear with time. \label{fig7}}
\end{figure}

Comparing the graphics from Figures~\ref{fig2}--\ref{fig7}, we can see that 
the asymptotic behavior is similar. However, the discrete-time model predicts, in some instances, 
smaller values. All the computations were done using \textsf{Mathematica}, version 12.1.

\vspace{-8pt}
\vspace{-8pt}
\vspace{-8pt}
\vspace{-8pt}
\vspace{-8pt}


\section{Discussion}
\label{sec4}

The COVID-19 pandemic has been a shocking experience to every human being, 
and no one expected the consequences felt worldwide. From this, we can also state that 
the number of variables and parameters to be taken into account in a \text{COVID-19} 
model are not small. Indeed, if there is such thing as a good model that explains COVID-19, 
surely it is not a simple model. Moreover, most of the models that explain epidemiological, 
ecological, and economical phenomena are not linear. Therefore, given the tools available, 
it is impossible to present an exact solution to the problem. That is one of the reasons 
for the development of nonstandard finite difference methods. Indeed, such methods 
enable us to construct discrete models from continuous ones, allowing one to present 
numerical solutions. 

Nowadays, there are several nonstandard finite schemes. 
We have used here the one developed by Mickens because it has shown to be dynamically consistent. 
Precisely, we discretized the continuous-time model presented and analyzed in \cite{PST}. 
With our discrete-time model, we analyzed the evolution of the COVID-19 pandemic in Portugal, 
from its beginning up to 4  May   2020. This choice allowed us to compare our results 
with those \mbox{of \cite{PST}}. We obtained the equilibrium points of the proposed discrete-time model,  
showing that they coincide with the ones from the continuous model. From the qualitative analysis 
point of view, we believe that it is important to prove the stability of the equilibrium points
if the model describes real data, which is our case here. For this reason, 
\mbox{Figures~\ref{fig2}--\ref{fig7}} show the simulation results made over 2000 days: our intention 
is to show the stability of the model. Regarding the stability, it should be noted that,
in contrast with \cite{PST}, which only proves local stability, here, we have proved 
the global stability of the disease-free equilibrium (DFE) and attempted to do the 
same with the endemic equilibrium (EE) point. However, regarding the EE, 
we were not able to prove the global stability analytically, 
and the question remains open. Figure~\ref{fig1} compares the predictions 
of the continuous model with those of the discrete-time model. 
One concludes that the discrete-time model fits a little bit better to the real data. 
From Figures~\ref{fig2}--\ref{fig7}, we present the convergence to the disease-free 
equilibrium point, using both continuous and discrete models. The asymptotic 
behavior is similar; the only difference worth mentioning is that the discrete-time model, 
in some small interval of time, predicts smaller values than the continuous one.

In this work, we restricted ourselves to the development of the pandemic 
in Portugal. We are aware that new models are emerging and, in a future work, 
we intend to analyze a model that fits the data of more than one country and region, 
possibly globally. This is under investigation and will be addressed elsewhere.
Another line of research concerns models that take vaccination into 
account \cite{Couras}. Here, vaccination is not considered because our model 
was created to explain the development of COVID-19 at the beginning of the pandemic. 
However, regarding vaccination, individuals do not become immune, 
so we think that we need to wait a little longer to see the development. 
It remains open the question of how to prove global stability 
for the endemic equilibrium. 


\authorcontributions{Conceptualization, S.V. and D.F.M.T.; 
methodology, S.V. and D.F.M.T.; software, S.V.;
validation, S.V. and D.F.M.T.; formal analysis, S.V. and D.F.M.T.; 
investigation, S.V. and D.F.M.T.; 
writing---original draft preparation, S.V. and D.F.M.T.; 
writing---review and editing, S.V. and D.F.M.T.; 
visualization, S.V. All authors have read and agreed 
to the published version of \mbox{the manuscript}.}

\funding{The authors were partially supported by 
the Portuguese Foundation for Science and Technology (FCT):
Sandra Vaz through the Center of Mathematics and Applications 
of \emph{Universidade da Beira Interior} (CMA-UBI), 
project UIDB/00212/2020; Delfim F. M. Torres through
the Center for Research and Development in Mathematics 
and Applications (CIDMA) of \emph{University of Aveiro}, 
project UIDB/04106/2020.}

\institutionalreview{Not applicable.}

\informedconsent{Not applicable.}

\dataavailability{The data and information used in this work are accessible 
to anyone and can be found in the following links: 
\url{https://www.pordata.pt/Portugal}  (accessed on 16 September 2021); 
\url{http://www.covid19.min-saude.pt} (accessed on 16 September 2021);   
\url{https://covid19.min-saude.pt/ponto-de situação atual-em-portugal} 
(accessed on 16 September 2021).} 

\acknowledgments{\textls[-15]{The authors are grateful 
to three reviewers for their several comments and suggestions.}}

\conflictsofinterest{The authors declare no conflict of interest.} 

\end{paracol}


\reftitle{References}


\end{document}